\documentclass[12pt]{article}
\usepackage{graphicx, tikz}

\usepackage{amssymb, amsmath}
\usepackage{amsthm}

\usepackage{comment}

\newcommand{\dd}{\mathrm{d}}

\newcommand{\p}{\mathbf{P}}

\theoremstyle{plain}
\newtheorem*{theorem}{Theorem}
\newtheorem{step}{Step}

\title{The Dynkin--Lamperti problem}

\author{P\'eter Kevei\footnote{Bolyai Institute, University of Szeged, Szeged,
Hungary. \\
E-mail: \texttt{kevei@math.u-szeged.hu}}}
\date{}

\begin{document}

\maketitle

\begin{abstract}
Let $X$ be a nonnegative random variable and $F(x) = \p ( X \leq x)$ its 
distribution function. It is well-known that if 
$\overline F(x) = 1 - F(x) = \ell(x) x^{-\alpha}$
as $x \to \infty$, where $\alpha \in (0,1)$ and $\ell$ is a slowly varying 
function, then $\lim_{x \to \infty} \overline F(x) U(x) = 
\frac{\sin \pi \alpha}{\pi \alpha}$, where $U(x) = \sum_{n=0}^\infty F^{*n} (x)$
is the renewal function of $F$. In this short note we prove the converse statement.
\end{abstract}

\section{Introduction and the result}

Let $X, X_1, X_2, \ldots $ be nonnegative
independent identically distributed (iid) random variables with distribution function $F (x ) = \p ( X \leq x) $. 
Let $S_n = \sum_{i=1}^n X_i$ denote the partial sum,  
$N_t = \max \{ k : S_k \leq t \}$ the corresponding 
renewal process, and $Y_t = t- S_{N_t}$ the spent lifetime,
and $Z_t = S_{N_t +1} - t$ the residual lifetime at $t$.
Let $U(x) = \sum_{n=0}^\infty F^{*n}(x)$ 
denote the renewal function, and put 
$\overline F (x) = 1 - F(x)$.

The Dynkin--Lamperti condition is
\begin{equation} \label{eq:dl1}
\overline F(x) = \frac{\ell(x)}{x^\alpha}, 
\end{equation}
where $\alpha \in (0,1)$, and $\ell$ is slowly varying.

In infinite mean renewal theory, the well-known Dynkin--Lamperti 
theorem (see, e.g.~Bingham et al.~\cite[Theorem 8.6.3]{BGT}) states that 
\eqref{eq:dl1} with $\alpha \in (0,1)$ is a  necessary and 
sufficient condition for the existence of the nondegenerate
distributional limit of any of $Y_t/t$, $Z_t/t$, $(Y_t/t, Z_t/t)$,
as $t \to \infty$. For the analogous result for subordinators,
see Kyprianou \cite[Section 5.5]{Kyprianou}.

Moreover, if \eqref{eq:dl1} holds then
\[
\lim_{x \to \infty} \overline F(x) \cdot U(x) = 
\frac{\sin (\pi \alpha)}{\pi \alpha}. 
\]
The Dynkin--Lamperti problem (see \cite[Section 8.6.3]{BGT}) is whether the converse implication holds. In this short note,
we prove that the converse indeed holds.

\begin{theorem}
Assume that  for some $c > 0$
\begin{equation} \label{eq:UF}
\lim_{x \to \infty} \overline F(x) \cdot U(x) = c. 
\end{equation}
Then $\overline F(x)$ is regularly varying 
with parameter $-\alpha$, for some 
$\alpha \in [0,1)$. If $\alpha > 0$ then it is the unique solution to
\begin{equation} \label{eq:alphac}
c = c_\alpha = 
\frac{\sin (\pi \alpha)}{\pi \alpha},
\end{equation}
while $c = c_0 = 1$ for $\alpha = 0$. In particular, $c \leq 1$.
\end{theorem}

The probabilistic interpretation
\[
\overline F(t) U(t) = 
\sum_{n=0}^\infty \p ( S_n \leq t, X_{n+1} > t) = 
\p ( Y_t + Z_t > t) 
\]
immediately shows that $c \leq 1$ in \eqref{eq:UF}.

Finally, we discuss the case $c = 0$.
Let
\begin{equation} \label{eq:defm}
m(x) = \int_0^x \overline F(u) \dd u
\end{equation}
denote the truncated mean. 
Then, monotonicity implies
\[
x = \int_0^x \overline F(x-y) U(y) \dd y \leq U(x) m(x).
\]
If \eqref{eq:UF} holds with $c = 0$, then for any $\lambda > 1$
\[
\begin{split}
\frac{m(\lambda x)}{m(x)} - 1 & =
\frac{\int_x^{\lambda x} \overline F(u) \dd u}{m(x)} 
\leq \frac{(\lambda -1) x \overline F(x)}{m(x)}  \\
& \leq 
(\lambda -1) \overline F(x) U(x) \to 0.
\end{split}
\]
That is, the truncated mean $m$ is slowly varying.


\section{Proof}

\begin{step}
Assume \eqref{eq:UF} and denote
\begin{equation} \label{eq:defh}
h(x) = \int_0^x \frac{\overline F(x-y)}{\overline F(y)} \dd y. 
\end{equation}
Then 
\begin{equation} \label{eq:hasy}
\lim_{x \to \infty} \frac{h(x)}{x} = \frac{1}{c}.
\end{equation}
\end{step}

\begin{proof}
For any $v \geq 0$
\[
\int_{[0,v]} \overline F(v-y) \dd U(y)  = 1.
\]
Integration on $(0,x)$ and Fubini's 
theorem give that
\[
\int_0^x \overline F(x-y) U(y) \dd y = x. 
\]
Since $U(y) \sim c / \overline F(y)$ as $y \to \infty$, we have
\[
\int_0^x \frac{\overline F(x-y)}{\overline F(y)} \dd y \sim \frac{x}{c},
\]
as claimed. Here, and later on, $\sim$ stands for 
asymptotic equality.
\end{proof}

The case $c=1$ is easier to handle.

\begin{step}
If $c = 1$ in \eqref{eq:UF} 
then $\overline F$ is slowly varying.
\end{step}

\begin{proof}
Changing variables we obtain 
\begin{equation} \label{eq:lambdah}
\lim_{x \to \infty} \int_0^1 
\frac{\overline F( x( 1 -u))}{\overline F(x u)} \dd u = c^{-1}.
\end{equation} 
Write
\[
\begin{split}
\int_0^1 \frac{\overline F( x( 1 -u))}{\overline F(x u)} \dd u 
& = \int_0^{1/2} \left(  
\frac{\overline F( x( 1 -u))}{\overline F(x u)} +
\frac{\overline F( xu)}{\overline F(x (1-u))} \right) \dd u.
\end{split}
\]
Since $a + a^{-1} \geq 2$ for any $a > 0$, this also 
implies that $c \leq 1$.
Let $c = 1$. We have $\overline F(xu) \geq \overline F(x(1-u))$, for $u \in (0,1/2)$, and $\frac{\overline F(xu)}{\overline F(x(1-u))}$
is nonincreasing on $(0,1/2)$.
Moreover, the function $y \mapsto y + y^{-1}$ is nondecreasing on 
$[1,\infty)$, therefore for $u_0 \in (0,1/2)$
\[
\begin{split}
& \int_0^{1/2} \left(  
\frac{\overline F( x( 1 -u))}{\overline F(x u)} +
\frac{\overline F( xu)}{\overline F(x (1-u))} -2  \right) \dd u \\
& \geq u_0 \left(  
\frac{\overline F( x( 1 -u_0))}{\overline F(x u_0)} +
\frac{\overline F( xu_0)}{\overline F(x (1-u_0)) } - 2 \right) \geq 0.
\end{split}
\]
As $x \to \infty$ the left-hand side tends to 0, thus  
for any $u_0 \in (0,1/2)$
\[
\lim_{x \to \infty} \frac{\overline F( x( 1 -u_0))}{\overline F(x u_0)} = 1,
\]
and the claim follows. 
\end{proof}

Therefore, the Theorem is proved for $c = 1$. 
In what follows, we assume that $c < 1$.

\begin{step} \label{step:Fbound}
There exist $c_1 > 0$, $0 < \beta < 1$, and $x_0 > 0$ such that 
for any $x \geq x_0$
\begin{equation} \label{eq:FFt}
\frac{\overline F(tx)}{\overline F(x)}
\begin{cases} 
\leq c_1 t^{-\beta}, & t \in (0,1], \\
\geq c_1^{-1} t^{-\beta}, & t > 1.
\end{cases}
\end{equation}
\end{step}

\begin{proof}
Recall \eqref{eq:defm}.
Then 
\begin{equation} \label{eq:hmF}
\begin{split}
h(2x) &  = \int_0^{2x} \frac{\overline F(2x- y)}{\overline F(y)} \dd y  
\geq \int_x^{2x} \frac{\overline F(2x- y)}{\overline F(y)} \dd y \\
&  \geq \frac{1}{\overline F(x)} \int_x^{2x} \overline F(2x- y) \dd y 
= \frac{m(x)}{\overline F(x)}.
\end{split}
\end{equation}
Furthermore,
\[
\overline F(x) \leq \frac{1}{x} \int_0^x \overline F(y) \dd y
= \frac{m(x)}{x}.
\]
Therefore, \eqref{eq:hmF} and \eqref{eq:hasy} imply 
that there exist $c_1 > 1$ and $x_0 > 0$ such that for $x \geq x_0$
\begin{equation} \label{eq:Fm}
\overline F(x) \leq \frac{m(x)}{x} \leq c_1 \overline F(x).
\end{equation}

Choose $x > y \geq x_0$. Then, noting that $m'(x) = \overline F(x)$ we have 
\[
\begin{split}
\log \frac{m(x)}{m(y)} & = \int_{y}^{x} \frac{m'(u)}{m(u)} \dd u \\
& \geq \frac{1}{c_1} \int_y^x u^{-1} \dd u = c_1^{-1} \log \frac{x}{y},
\end{split}
\]
implying that 
\begin{equation} \label{eq:FF}
\frac{\overline F(y)}{\overline F(x)} \leq 
c_1 \frac{m(y)}{y} \frac{x}{m(x)} 
 \leq c_1 \left( \frac{y}{x} \right)^{c_1^{-1} -1}.
\end{equation}
Put $\beta = 1 - c_1^{-1} > 0$.
Therefore, for $t \in (0,1]$, if $x \geq x_0$ and $tx \geq x_0$ then
\[
\frac{\overline F(tx)}{\overline F(x) } \leq c_1 t^{- \beta},
\]
while if $tx \leq x_0$, then
\[
\frac{\overline F(tx)}{\overline F(x) } 
\leq \frac{1}{\overline F(x_0)}
\frac{\overline F(x_0)}{\overline F(x)} 
\leq \frac{c_1}{\overline F(x_0)}
\left( \frac{x_0}{x} \right)^{-\beta} 
\leq \frac{c_1}{\overline F(x_0)} t^{-\beta}.
\]
Therefore, increasing $c_1$ if necessary, the statement follows.

For $t > 1$, the bound follows by substituting $x' = tx$ and 
$t' = t^{-1}$ in the first bound.
\end{proof}

This step allows us to use the selection principle. For any 
sequence $x_n \to \infty$ we can choose a subsequence (still denoted by $x_n$)
and there exists a nondecreasing function $V$ such that 
\begin{equation} \label{eq:defV}
\lim_{n \to \infty} \frac{\overline F(x_n)}{\overline F(tx_n)} = V(t),
\end{equation}
for each $t \in (0,\infty)$, which is a continuity point of $V$.
Our aim is to show that $V(t) =t^\alpha$ for $\alpha \in (0,1)$,
determined via \eqref{eq:alphac}, is the unique subsequential limit.
First, we show some basic properties of $V$.

\begin{step} \label{step:Vreneq}
Let $V$ be any subsequential limit in \eqref{eq:defV}. Then 
$\lim_{t \downarrow 0}V(t) = 0$,  $V(1) = 1$, $V$ is continuous and nondecreasing,
and it satisfies the equation
\begin{equation} \label{eq:Vreneq}
\int_0^t \frac{V(u)}{V(t-u)} \dd u = \frac{t}{c}, 
\quad t \in (0,\infty).
\end{equation}
Furthermore, $V$ satisfies the renewal equation 
\begin{equation} \label{eq:Vreneq2}
\int_{[0,t]} \frac{1}{V(t-s)} \dd V(s) = c^{-1},
\end{equation}
for every $t \in (0,\infty)$.
\end{step}

\begin{proof}
Writing $t x_n$ in \eqref{eq:hasy} and changing variables, we see that
\[
\lim_{n \to \infty} \frac{1}{t} 
\int_0^t \frac{\overline F(x_n(t-u))}{\overline F(x_n u)}\dd u
= c^{-1}.
\]
The integrand converges pointwise almost everywhere 
to $\frac{V(u)}{V(t-u)}$, and 
the estimate \eqref{eq:FFt} allows us to use Lebesgue's dominated 
convergence theorem; thus, \eqref{eq:Vreneq} follows. From \eqref{eq:defV}, 
it is immediate that $V$ is nondecreasing; thus, $\lim_{t \downarrow 0} V(t) = V(0)$ exists. If $V(0) > 0$, then 
\[
\lim_{t \downarrow 0} t^{-1} \int_{0}^t \frac{V(u)}{V(t-u)} \dd u = 1,
\]
which contradicts the fact that $c < 1$.
Thus, $V(0) = 0$.

Fubini's theorem gives that 
\[
\int_{0}^t \int_{[0,s]} \frac{1}{V(s-y)} \dd V(y) \dd s 
= \int_0^t \frac{V(s)}{V(t-s)} \dd s,
\]
and the latter equals $\frac{t}{c}$ by \eqref{eq:Vreneq}. 
Thus,  \eqref{eq:Vreneq2} holds for almost every $t > 0$. 
For any $a > 0$
\[
\int_{[0,s]} \frac{1}{V(s-y)} \dd V(y) \geq \frac{V(a+)-V(a-)}{V(s-a)},
\]
and the latter tends to $\infty$ as $s \downarrow a$ if 
$V(a+) - V(a-) > 0$. Therefore, $V$ is continuous. 
From \eqref{eq:defV} together with the continuity, we see that $V(1) = 1$.

Finally, we show that \eqref{eq:Vreneq2} holds for 
every $t \in (0,\infty)$. Put
\[
\mu (t) = \int_{[0,t]} \frac{1}{V(t-s)} \dd V(s).
\]
Since $\mu(t) = c^{-1}$ for almost every $t \in (0,\infty)$, for 
any $t \in (0,\infty)$ we can choose a subsequence $t_n \uparrow t$,
such that $\mu(t_n) = c^{-1}$. Thus, by Fatou's lemma 
\begin{equation} \label{eq:mu-bound}
\mu(t) \leq c^{-1} \quad \text{for every } \ t \in (0,\infty).
\end{equation}
Next, put 
\begin{equation} \label{eq:def-mV}
m_V(t) = \int_0^t \frac{1}{V(s)} \dd s.
\end{equation}
Monotonicity gives that
\[
\frac{2t}{c} = \int_0^{2t} \frac{V(u)}{V(2t-u)} \dd u \geq 
V(t) m_V(t).
\]
This implies that 
\begin{equation} \label{eq:mV-bound}
m_V(t) \leq \frac{2t}{c V(t)}, \text{ for all } \ t \in (0,\infty).
\end{equation}
For $t \in (0,\infty)$, $h > 0$, let
\begin{equation} \label{eq:def_A}
A_h(t) = \frac{1}{h} \int_{(t-h)_+}^t \frac{1}{V(u)} \dd u.
\end{equation}
For $h \geq t$, 
\[
A_h(t) = \frac{1}{h} \int_0^t \frac{1}{V(u)} \dd u \leq
\frac{m_V(t)}{t},
\]
while for $h \leq t$, the monotonicity of $V$ implies
the same upper bound. Combining with \eqref{eq:mV-bound}, we 
obtain that  
\begin{equation} \label{eq:A-bound}
A_h(t) \leq \frac{2}{c V(t)}.
\end{equation}
Using Fubini's theorem
\[
\frac{1}{c} = \frac{1}{h} \int_{t-h}^t  \mu(u) \dd u 
= \int_{[0,t]} A_h(t-s) \dd V(s).
\]
Fix $t > 0$. Since $V$ is continuous, 
$\lim_{h \downarrow 0} A_h(t-s)  = V(t-s)^{-1}$ for 
all $s \in (0,t)$. 
Using the bound \eqref{eq:A-bound}, Lebesgue's dominated 
convergence theorem gives that \eqref{eq:Vreneq2} holds 
for all $t \in (0,\infty)$.
\end{proof}

The following comparison argument is the key step towards 
the proof of uniqueness.

\begin{step} \label{step:comparison}
Let $V$ and $W$ be continuous, nondecreasing functions, 
$V(0) =W(0) = 0$, 
that satisfy the renewal equation \eqref{eq:Vreneq2} for every $t \in (0,\delta]$ for some $\delta > 0$. Assume that 
$V(s) \leq W(s)$ for $s \in [0,\delta]$ and 
$V(\delta ) = W(\delta)$. Then $V(s) = W(s)$ for $s \in [0,\delta]$.
\end{step}

\begin{proof}
Let 
\[
Z(t) = (V * W)(t) =  \int_{[0,t]} V(t-s) \dd W(s).
\]
Then
\[
\left(\frac{1}{V} * V * W \right)(t) = c^{-1} W(t),
\]
and similarly
\[
\left(\frac{1}{W} * W * V \right)(t) = c^{-1} V(t).
\]
Since convolution is commutative, taking the difference
\begin{equation} \label{eq:VW}
\int_{[0,t]} \left( \frac{1}{V(t-s)} - \frac{1}{W(t-s)} \right)
\dd Z(s) = c^{-1} (W(t) - V(t)).
\end{equation}
Since $V(s) \leq W(s)$ for all $s \in (0,\delta)$ and 
$V(\delta) = W(\delta)$,
\[
\frac{1}{V(\delta -s)} \geq  \frac{1}{W(\delta -s)}
\]
for all $s \in (0,\delta)$, and its integral on $[0,\delta]$
with respect to $Z$ is 0. Therefore, 
$V(\delta -s) = W(\delta -s)$ for each $s$ which is a point 
of increase of $Z$. By \eqref{eq:Vreneq2}, 0 is a point of increase
of $V$ and $W$, and thus also of $Z$. Therefore, 
$V(\delta -s_0) = W(\delta -s_0)$ for some $s_0 \in (0,\delta)$.
Repeating this argument, we get that $V(\delta - ks_0) = W(\delta - k s_0)$ for $k\in \{1,2,\ldots \}$, with $\delta - k s_0 \geq 0$. Using that $s_0$ can be chosen arbitrarily 
small, and that $V$ and $W$ are continuous, we obtain that 
$V(s) = W(s)$ for all $s \in [0,\delta]$, as claimed.
\end{proof}

\begin{step} \label{step:sandwich}
Let $V$ be a continuous, nondecreasing function, $V(0) = 0$ and $V(1) = 1$,
that satisfies the renewal equation \eqref{eq:Vreneq2}
for $t \in (0,1]$. Assume further, that for some $\rho \in (0,1)$ there 
exist $0< C_1 < C_2 < \infty$ such that
\[
C_1 t^\rho \leq V(t) \leq C_2 t^{\rho}, \quad t \in (0,1].
\]
Then $\rho = \alpha$, where $\alpha \in (0,1)$
satisfies $c = \frac{\sin (\pi \alpha)}{ \pi \alpha}$.
\end{step}

\begin{proof}
Choose $\alpha \in (0,1)$ such that 
$c = c_\alpha = \frac{\sin(\pi \alpha)}{\pi \alpha}$.
Then $W(t) = t^\alpha$ satisfies \eqref{eq:Vreneq2}.

Assume that $\rho < \alpha$. Then 
$\lim_{t \to 0} \frac{V(t)}{t^\alpha} = \infty$; thus, 
there exists $\delta \in (0,1]$ such that $V(t) \geq t^\alpha$ for 
$t \in (0,\delta)$ and $V(\delta ) = \delta^\alpha$. Then
Step \ref{step:comparison} implies that $V(s) = s^\alpha$ for 
$s \in (0,\delta]$, contradicting $\rho < \alpha$.

Similarly, if $\rho > \alpha$, then 
$\lim_{t \to 0} \frac{V(t)}{t^\alpha} = 0$, and the same comparison 
argument leads to a contradiction. Therefore $\rho = \alpha$.
\end{proof}

\begin{step} 
Let $V$ be a continuous, nondecreasing function, $V(0) = 0$ and $V(1) = 1$,
that satisfies the renewal equation \eqref{eq:Vreneq2}
for $t \in (0,\infty)$. Then $V(t) = t^{\alpha}$, $t \in [0,1]$,
where $\alpha \in (0,1)$
satisfies 
$c = \frac{\sin (\alpha \pi)}{\alpha \pi}$.
\end{step}

\begin{proof}
Fix $\lambda \in (0,1)$, and let
\begin{equation} \label{eq:defL}
L = \sup_{0 < u \leq 1} \frac{V(\lambda u)}{V(u)}.
\end{equation}
Choose a subsequence $u_n \in (0,1)$ such that 
\begin{equation} \label{eq:defu}
\lim_{n \to \infty} u_n = u \in [0,1], \quad \text{and } \ 
\lim_{n \to \infty} \frac{V(\lambda u_n)}{V(u_n)} = L.
\end{equation}
Let 
\begin{equation} \label{eq:defG}
G_n(t) = \frac{V(u_n t)}{V(u_n)}, \quad t \in (0,\infty).
\end{equation}

We use the same argument as in Step \ref{step:Fbound}
to bound $G_n$. 
Monotonicity of $V$ implies that 
$m_V(t) \geq \frac{t}{V(t)}$; therefore, using
also \eqref{eq:mV-bound}
\begin{equation} \label{eq:Vbound}
\frac{t}{m_V(t)} \leq V(t) \leq \frac{2}{c} \frac{t}{m_V(t)}.
\end{equation}
Next, for $0 < a < b < \infty$
\[
\log \frac{m_V(b)}{m_V(a)} = 
\int_a^b \frac{m_V'(t)}{m_V(t)} \dd t 
\geq \frac{c}{2} \log \frac{b}{a},
\]
implying that
\[
\frac{m_V(b)}{m_V(a)} \geq \left( \frac{b}{a} \right)^{\frac{c}{2}}.
\]
Combining the latter inequality with \eqref{eq:Vbound},
we obtain for $0 < y < x < \infty$
\begin{equation} \label{eq:Vratiobound}
1 \leq \frac{V(x)}{V(y)} \leq \frac{2x}{c m_V(x)} \frac{m_V(y)}{y}
\leq \frac{2}{c} \left( \frac{x}{y} \right)^{1 - \frac{c}{2}}.
\end{equation}
Thus 
\begin{equation} \label{eq:Gbound}
\begin{split}
G_n(t) & = \frac{V(u_n t)}{V(u_n)} \geq \frac{c}{2} t^{1- \frac{c}{2}}, \quad t \in (0,1], \\
G_n(t) & = \frac{V(u_n t)}{V(u_n)} \leq \frac{2}{c} t^{1- \frac{c}{2}}, \quad t \in [1,\infty), 
\end{split}
\end{equation}

There exists a subsequence $n_k$ such that 
\begin{equation} \label{eq:defG-limit}
\lim_{k \to \infty} G_{n_k} (t) = G(t),
\end{equation}
for all continuity points of $G$. Since
\[
\int_0^t \frac{G_n(s)}{G_n(t-s)} \dd s = 
\int_0^t \frac{V(u_n s)}{V(u_n (t-s))} \dd s = \frac{t}{c},
\]
and \eqref{eq:Gbound} gives an integrable bound,  
Lebesgue's dominated convergence theorem implies that 
\[
\int_{0}^t \frac{G(s)}{G(t-s)} \dd s = \frac{t}{c}, 
\quad t \in (0,\infty).
\]
By Step \ref{step:Vreneq}, this implies that 
$G$ is continuous. Furthermore, the choice of 
$(u_n)$ in \eqref{eq:defu} shows that 
$G(\lambda) = L$.

Let 
\begin{equation} \label{eq:tildeG}
\widetilde G (t) = \frac{G(\lambda t)}{L}, \quad t \in [0,1].
\end{equation}
Then $\widetilde G$ is a solution to \eqref{eq:Vreneq2}, and 
$\widetilde G(1) = 1$. Furthermore, by the choice of $u_n$
\[
\widetilde G(t) \leq G(t), \quad t \in (0,1],
\]
which implies by Step \ref{step:comparison} that 
$\widetilde G(t) = G(t)$, that is
\begin{equation} \label{eq:Gscale}
G(\lambda t) = L G(t), \quad t \in [0,1].
\end{equation}
Iteration gives that $G(\lambda^k) = L^k$, thus 
with $\rho = \frac{\log L}{\log \lambda}$
\[
L^2 t^{\rho} \leq G(t) \leq L^{-1} t^\rho, \quad t \in (0,1].
\]
It follows from \eqref{eq:Gbound} that $\rho \in (0,1)$.
Step \ref{step:sandwich} implies $\rho = \alpha$, 
where $\alpha \in (0,1)$ is the unique solution to $c_\alpha = c$. Thus 
$L = \lambda^\alpha$, and we obtained that 
\begin{equation} \label{eq:Vsup}
\sup_{0 < u \leq 1} \frac{V(\lambda u)}{V(u)} = \lambda^\alpha.
\end{equation}
The same argument can be repeated to show that 
\begin{equation} \label{eq:Vinf}
\inf_{0 < u \leq 1} \frac{V(\lambda u)}{V(u)} = \lambda^\alpha.
\end{equation}

These bounds together show that $\frac{V(\lambda u)}{V(u)} = \lambda^\alpha$. Since this holds for any $\lambda \in (0,1)$, choosing $u=1$, the statement follows.
\end{proof}

A similar argument shows that $V(t) = t^\alpha$ holds for any
$t > 0$, but for our purposes, the uniqueness on $[0,1]$ is 
enough.

Therefore, the limit in \eqref{eq:defV} is unique, which 
implies that 
\[
\lim_{x \to \infty} \frac{\overline F(x)}{\overline F(tx)} =
t^{\alpha},
\]
which proves the theorem.

\bigskip 

\noindent \textbf{AI declaration.}
Around 2021, I tried to solve the problem, and Step 1 seemed to be a promising start. The case $c = 1$ is due to Professor Vilmos Totik. 
All the subsequent
steps, ideas, and proofs are based on my conversations with 
Chat-GPT 6.1. 
I take full responsibility for the content of the paper,
the correctness of its mathematical results and proofs, and any remaining errors.


\begin{thebibliography}{1}

\bibitem{BGT}
N.~H. Bingham, C.~M. Goldie, and J.~L. Teugels.
\newblock {\em Regular Variation}, volume~27 of {\em Encyclopedia of Mathematics and its Applications}.
\newblock Cambridge University Press, Cambridge, 1989.

\bibitem{Kyprianou}
A.~E. Kyprianou.
\newblock {\em Fluctuations of {L}\'{e}vy processes with applications}.
\newblock Universitext. Springer, Heidelberg, second edition, 2014.

\end{thebibliography}

\end{document}